\documentclass[11pt]{amsart}
\usepackage[T1]{fontenc}
\usepackage{lmodern}
\usepackage{amsmath,amssymb,amsthm,mathtools}
\usepackage{microtype}
\usepackage{enumitem}
\usepackage[colorlinks=true,linkcolor=blue,citecolor=blue,urlcolor=blue]{hyperref}
\newtheorem{theorem}{Theorem}[section]
\newtheorem{lemma}[theorem]{Lemma}
\newtheorem{proposition}[theorem]{Proposition}

\theoremstyle{definition}

\theoremstyle{remark}
\newtheorem{remark}[theorem]{Remark}
\newcommand{\F}{\mathbb F}
\newcommand{\Z}{\mathbb Z}

\newcommand{\PSL}{\operatorname{PSL}}
\newcommand{\PSp}{\operatorname{PSp}}
\newcommand{\SL}{\operatorname{SL}}
\newcommand{\Aut}{\operatorname{Aut}}
\newcommand{\Sym}{\operatorname{Sym}}

\newcommand{\ord}{\operatorname{ord}}
\newcommand{\tr}{\operatorname{tr}}

\title[Four Kourovka Notebook problems]{On four problems from the Kourovka Notebook:\\ spectra, periodic automorphism groups, commutators, and operator simplicity}
\author{Ivan Samuil Nnuga Mongo}
    
\address{
Novosibirsk State University,
Novosibirsk, Russia
}
\email{i.nnuga.mongo@g.nsu.ru}

\begin{document}
\begin{abstract}
We examine Problems 18.62(a), 18.63, 18.65, and 20.120 of the
Kourovka Notebook.  First, a family of semidirect products involving
$\PSL_2(p^2)$ gives arbitrarily many pairwise nonisomorphic groups
critical with respect to one simple-group spectrum, subject to explicit
arithmetic congruences.  Second, a completed-algebra construction
using sparse torsion identities and a matrix Golod--Shafarevich--Vinberg
inequality produces an infinite periodic group generated by two
fixed-point-free automorphisms of order five of an abelian group.
Third, we give a short matrix proof that the class of elements of
order three in $\PSL_2(8)$ contradicts the assertion in Problem 18.65
as written.  Finally, we spell out how P. Hall's 1962 wreath-power
theorems answer Problem 20.120 under the usual normal-subgroup meaning
of simplicity for a group with one distinguished automorphism.
The four arguments have substantially different provenance.  In
particular, the last is a consequence of a known construction, not a
priority claim.  The two longer constructions are recorded in detail
for independent mathematical scrutiny.
\end{abstract}
\maketitle

\section{Statements, conventions, and limits of the claims}
For a finite group $X$ let $\omega(X)$ be the set of orders of its
elements.  A \emph{section} of $X$ is $H/N$, where $N\trianglelefteq
H\le X$.  A group is $\omega$-\emph{critical} if its spectrum is
$\omega$ and no proper section has this spectrum.  We count groups up
to isomorphism.  The four problems considered here are independent:
Problem 18.62(a) asks for a universal bound on the number of
$\omega(L)$-critical groups as $L$ ranges over finite simple groups;
18.63 concerns two fixed-point-free order-five automorphisms of an
abelian group generating a periodic group; 18.65 asserts a normal
$p$-radical conclusion for a class of $p$-elements whose pairwise
commutators are $p$-elements, with $p\ne5$; and 20.120 concerns
operator simplicity of a locally nilpotent group \cite{Notebook}.
No claim is made here to solve the classification problem 18.62(b).

\begin{theorem}[18.62(a)]\label{thm:1862}
For every positive integer $N$ there exist an odd prime $p$ and $N$
pairwise nonisomorphic $\omega(\PSp_4(p))$-critical finite groups.
Consequently there is no universal upper bound as in 18.62(a).
\end{theorem}
\begin{theorem}[18.63]\label{thm:1863}
There exist an abelian group $A$ and fixed-point-free automorphisms
$\alpha,\beta\in\Aut(A)$, each of order five, such that
$\langle\alpha,\beta\rangle$ is infinite and periodic.
\end{theorem}
\begin{theorem}[18.65]\label{thm:1865}
Let $G=\PSL_2(8)$, $p=3$, and $C$ be its conjugacy class of elements
of order three.  Every commutator $[c,d]$, $c,d\in C$, has order one
or nine, whereas $O_3(G)=1$.  Hence the assertion of 18.65, as
stated with ``$p$-element'', is false.
\end{theorem}
\begin{theorem}[20.120; consequence of Hall]\label{thm:20120}
For any prime $p$ there is a countable, locally finite $p$-group $D$
and an automorphism $\sigma$ of infinite order such that $D$ is
perfect, every normal $\sigma$-invariant subgroup of $D$ is $1$ or
$D$, and $C_D(\sigma^m)=1$ for every nonzero integer $m$.
\end{theorem}

\section{Critical groups for the spectrum of
\texorpdfstring{$\PSp_4(p)$}{PSp4(p)}}
\subsection{The arithmetic choice of prime}
Fix $N\ge1$ and set $B=2N+1$.  Let $\Pi_B$ contain every odd prime
factor of $s^2+1$ for odd $1\le s\le B$, as well as every odd prime
factor of $s^2-1$ for odd $3\le s\le B$.  Set
\[
 M_B=4\prod_{\ell\in\Pi_B}\ell.
\]
Dirichlet's theorem yields a prime $p>B+2$ with $p\equiv1\pmod{M_B}$.
The following elementary estimate is the reason for this choice.
\begin{lemma}\label{lem:arith}
For $s\in\{-B,-B+2,\ldots,B\}$ and $\varepsilon\in\{\pm1\}$,
put $u=\varepsilon+ps$.  Then
\[
 \gcd(p^2+1,u)=2,
\]
and the integer $\gcd(p^2-1,u)/2$ divides either $p-1$ or $p+1$.
\end{lemma}
\begin{proof}
If an odd prime $\ell$ divides $p^2+1$ and $u$, then
$p^2\equiv-1$ and $ps\equiv-\varepsilon\pmod\ell$, whence
$\ell\mid s^2+1$.  Our choice gives $p\equiv1\pmod\ell$, a
contradiction.  Now $v_2(p^2+1)=1$ and $u$ is even, proving the
first assertion.

If an odd prime $\ell$ divides $p^2-1$ and $u$, then
$p\equiv\pm1\pmod\ell$.  The case $p\equiv-1$ implies
$s\equiv\varepsilon\pmod\ell$.  Unless $s=\varepsilon$, this forces
$\ell\mid s^2-1$ (the case $s=-\varepsilon$ produces only the factor
two), contradicting $p\equiv1\pmod\ell$.  Thus, when
$s\ne\varepsilon$, the odd part of the gcd divides $p-1$.
Since $p\equiv1\pmod4$, $v_2(p+1)=1$; after dividing the gcd
by two, its remaining $2$-part also divides $p-1$.
When $s=\varepsilon$, $u=\varepsilon(p+1)$ and the gcd is $p+1$.
\end{proof}

\subsection{Modules and element orders}
Let
\[
 E=\F_{p^2},\qquad S=\PSL_2(E),\qquad
 H=S\rtimes\langle\tau\rangle,
\]
where $\tau$ is the field automorphism $x\mapsto x^p$ of order two.
For any odd $b\in\{3,5,\ldots,B\}$ put
\[
 W_b=U\otimes_E \Sym^b(U)^{(p)},
\]
with $U$ the natural two-dimensional module for $\SL_2(E)$.
Here the superscript $(p)$ denotes precomposition with the
entrywise $p$-power automorphism.  Since $b<p$ and $b$ is odd,
Steinberg's tensor-product classification
\cite[Lemma 6]{Lytkin2018} identifies $W_b$ as a simple
$ES$-module: the central involution of $\SL_2(E)$ acts by
$(-1)^{1+b}=1$.  Set $V_b=\operatorname{Res}_{E/\F_p}W_b$.
Entrywise Frobenius on the usual tensor basis defines a
$\F_p$-linear involution $F$ satisfying
$F\rho(g)F^{-1}=\rho(\tau(g))$, hence extends $V_b$ to an
$\F_pH$-module.  Over an algebraic closure,
\[
 V_b\otimes_{\F_p}\overline{\F}_p
 \cong W(1,b)\oplus W(b,1).
\]
The two constituents are nonisomorphic because $b>1$, and the
outer involution interchanges them.  It follows that $V_b$ is
irreducible as an $\F_pH$-module and
\begin{equation}\label{eq:dimension}
 \dim_{\F_p}V_b=4(b+1).
\end{equation}
The representation is faithful: the nontrivial simple group $S$
acts nontrivially, and $\tau$ interchanges its two constituents.

We use the following standard element-order formulas, recorded in
\cite[Lemmas 3 and 5]{Lytkin2018}, for $p>3$:
\begin{equation}\label{eq:specpsp}
 \omega(\PSp_4(p))=\operatorname{Div}\left\{
 p(p-1),\ p(p+1),\ \frac{p^2-1}{2},\ \frac{p^2+1}{2}\right\},
\end{equation}
where $\operatorname{Div}$ denotes the set of positive divisors,
and
\begin{equation}\label{eq:spech}
 \omega(H)=\omega(S)\cup 2\omega(\PSL_2(p))
 \subseteq\omega(\PSp_4(p)).
\end{equation}

\begin{lemma}[Orders in the field-automorphism coset]
\label{lem:fieldcoset}
Every $p'$-element of $H\setminus S$ has order dividing
$p-1$ or $p+1$, and both orders $p-1$ and $p+1$ occur
in this coset.  Every element of this coset whose order
is divisible by $p$ has order $2p$.
\end{lemma}
\begin{proof}
Write $h=\overline g\tau$ with $g\in\SL_2(E)$, and put
$n=gg^{(p)}$, where $(p)$ denotes entrywise Frobenius.
Then $h^2=\overline n$ and
\[
 \tr(n)^p=\tr(g^{(p)}g)=\tr(gg^{(p)})=\tr(n).
\]
Thus $\tr(n)\in\F_p$ and $\det(n)=1$.
If $n$ is semisimple, its eigenvalues either belong to
$\F_p^\times$ or form a Frobenius-conjugate pair of norm
one in $E^\times$.  Hence $r=\ord(n)$ divides $p-1$ or
$p+1$.  Since $h$ is outside the subgroup of index two,
\[
 \ord(h)=2\ord(\overline n)
        =\frac{2r}{\gcd(r,2)}.
\]
If $r$ is even, this equals $r$; if $r$ is odd, $2r$
still divides the corresponding even integer $p-1$ or
$p+1$.  If $n$ is not semisimple, its projective image
is a nonidentity unipotent element of order $p$,
and therefore $\ord(h)=2p$.

To realize $p-1$, choose $\lambda\in\F_p^\times$ of order
$p-1$ and $x\in E^\times$ with $xx^p=\lambda$.
Such $x$ exists by surjectivity of the norm map.
For $g=\operatorname{diag}(x,x^{-1})$ we have
$gg^{(p)}=\operatorname{diag}(\lambda,\lambda^{-1})$,
so $\overline g\tau$ has order $p-1$.

To realize $p+1$, choose $\lambda\in E^\times$ of order
$p+1$ and put $t=\lambda+\lambda^{-1}\in\F_p$.
Since $\lambda\ne-1$, we have $t+2\ne0$.
Choose $x\in E^\times$ with $xx^p=t+2$ and set
\[
 g=\begin{pmatrix}0&-1\\1&x\end{pmatrix}.
\]
Then $\det(g)=1$ and $\tr(gg^{(p)})=xx^p-2=t$.
Thus $gg^{(p)}$ has eigenvalues $\lambda,\lambda^{-1}$
and order $p+1$, so $\overline g\tau$ has order $p+1$.
\end{proof}

\begin{lemma}\label{lem:isospectral}
For every indicated $b$, the semidirect product
$G_b=V_b\rtimes H$ has the same spectrum as $\PSp_4(p)$.
\end{lemma}
\begin{proof}
For a semisimple lift $\widetilde h\in\SL_2(E)$ with
eigenvalues $\lambda,\lambda^{-1}$, the weights on $W(1,b)$ are
\begin{equation}\label{eq:weights}
 \lambda^{\varepsilon+ps},\qquad
 \varepsilon\in\{\pm1\},\quad s=-b,-b+2,\ldots,b.
\end{equation}
On $W(b,1)$ they are $\lambda^{s+p\varepsilon}$.
Because $p^2\equiv-1\pmod{p^2+1}$, multiplication by $p$
transforms the latter exponents into exponents of the former shape
up to signs.  Lemma~\ref{lem:arith} shows that a nonidentity
semisimple element of a nonsplit torus of $S$ has no fixed vector
on $V_b$; for a semisimple $h\in S$ in a split torus,
\begin{equation}\label{eq:fixedrestriction}
 C_{V_b}(h)\ne0\quad\Longrightarrow\quad
 \ord(h)\mid p-1\ \text{or}\ \ord(h)\mid p+1.
\end{equation}
Indeed, eigenvalues in a split (respectively nonsplit) torus have
orders dividing $p^2-1$ (respectively $p^2+1$), while the image in
$S$ identifies the central signs.

For a unipotent $u\in S$ of order $p$, the nilpotence index of
$\rho(u)-1$ on $W_b$ is at most $b+2<p$.  Consequently
\[
 1+\rho(u)+\cdots+\rho(u)^{p-1}
 =(\rho(u)-1)^{p-1}=0.
\]
By Lemma~\ref{lem:fieldcoset}, elements of $H\setminus S$
with nontrivial $p$-part have order $2p$; their square has
order $p$, and the same geometric sum
argument gives $1+h+\cdots+h^{2p-1}=0$ as an endomorphism of
$V_b$.  Thus the extension introduces no elements of order $p^2$
or $2p^2$.

For a $p'$-element $h\in H$ of order $m$, the only possible
orders above $h$ in $V_b\rtimes H$ are $m$ and $pm$, and the
second is possible precisely when $C_{V_b}(h)\ne0$.
If $h\in S$, equation \eqref{eq:fixedrestriction} puts $pm$ in
\eqref{eq:specpsp}; if $h\notin S$, Lemma~\ref{lem:fieldcoset}
shows that $m\mid p-1$ or $m\mid p+1$.
We have proved $\omega(G_b)\subseteq\omega(\PSp_4(p))$.

The maximal orders $(p^2\pm1)/2$ are already represented in $S$.
By Lemma~\ref{lem:fieldcoset}, choose $t_-,t_+\in H\setminus S$
of orders $p-1$ and $p+1$, respectively.
Their squares are in $S$, of orders $(p-1)/2$ and $(p+1)/2$.
For $t_-^2$ one may take an eigenvalue $\lambda$ of order
$p-1$ with $\lambda^p=\lambda$; the weight
$1-p$ in \eqref{eq:weights} is then trivial.  For $t_+^2$ one
may take $\lambda$ of order $p+1$, with
$\lambda^p=\lambda^{-1}$; the weight $1+p$ is trivial.
The involution $\tau$ exchanges the two absolutely simple
constituents.  Thus from a nonzero vector $w$ fixed by $t_\pm^2$
in one constituent, $w+t_\pm w$ is a nonzero vector fixed by
$t_\pm$.  Therefore the semidirect product contains elements
of orders $p(p-1)$ and $p(p+1)$.  All four maximal element
orders in \eqref{eq:specpsp} occur, proving equality.
\end{proof}

\subsection{Criticality}
\begin{lemma}\label{lem:critical}
Each $G_b$ is $\omega(\PSp_4(p))$-critical.
\end{lemma}
\begin{proof}
Write $V=V_b$ and $G=V\rtimes H$.  Suppose a subgroup $M\le G$
has the full spectrum.  Its projection to $H$ contains elements
of orders
\[
 A=\frac{p^2-1}{2},\qquad B'=\frac{p^2+1}{2}.
\]
Both elements project to $S$, because they are $p'$-elements
and Lemma~\ref{lem:fieldcoset} excludes both orders in the
outer coset: both $A$ and $B'$ exceed $p+1$.
By the classification of maximal subgroups of $\PSL_2(p^2)$
(Dickson's theorem \cite[Chapter XII]{Dickson1901}), a proper
subgroup cannot contain both
torus orders $A$ and $B'$ when $p>5$: a Borel subgroup or a
split torus normalizer misses $B'$, a nonsplit torus normalizer
misses $A$, a subfield group misses $B'$, and an exceptional
subgroup $A_4$, $S_4$, or $A_5$ has insufficient order.
Consequently $S\le\pi(M)$, and either $\pi(M)=S$ or
$\pi(M)=H$.

If $\pi(M)=S$, then $M\le V\rtimes S$.
But this group does not have elements of order $p(p-1)$ or
$p(p+1)$.  Indeed a semisimple element of $S$ of order $p-1$
(resp. $p+1$) lifts to eigenvalues in $\SL_2(E)$ of order
$2(p-1)$ (resp. $2(p+1)$).  To fix a weight in
\eqref{eq:weights}, the exponent $\varepsilon+ps$ would have
to be divisible by the corresponding order.  Reducing modulo
$p\mp1$ gives $s=\mp\varepsilon$ since $|s|\le B<p-2$;
the exponent is then only $\pm(p\mp1)$ and fails the
required divisibility by $2(p\mp1)$.
The same holds for the Galois-conjugate constituent.
Thus $\pi(M)\ne S$.

If $\pi(M)=H$, then $M\cap V$ is $H$-stable, hence either $0$
or $V$ by irreducibility.  In the first case $M\cong H$ and
cannot have the full spectrum; in the second $M=G$.
Hence no proper subgroup has the full spectrum.

Finally, since $V$ is faithful and irreducible, every nontrivial
normal subgroup $N$ of $G$ contains $V$: if $N\cap V=1$,
then $N$ centralizes $V$ and lies in $C_G(V)=V$, a contradiction.
The quotient $G/V\cong H$ does not have the full spectrum.
Therefore no nontrivial quotient does either.  Since a section
$K/N$ has its spectrum contained in $\omega(K)$, the two facts
prove criticality.
\end{proof}
\begin{proof}[Proof of Theorem~\ref{thm:1862}]
For each $N$, choose $p$ as in Lemma~\ref{lem:arith}, and take
the $N$ modules indexed by $b=3,5,\ldots,2N+1$.
Lemmas~\ref{lem:isospectral} and \ref{lem:critical} give the desired
critical groups.  They are pairwise nonisomorphic because
$\dim_{\F_p}V_b=4(b+1)$ are distinct, and $V_b$ is the unique
minimal nontrivial normal subgroup of $G_b$.
\end{proof}

\section{Two fixed-point-free automorphisms of order five}
\subsection{Cyclotomic algebra and its formal neighbourhood}
Put $K=\F_{11}$ and let
\begin{equation}\label{eq:freeproduct}
 U=K\langle a,b\rangle/(\Phi_5(a),\Phi_5(b)),
 \qquad \Phi_5(X)=1+X+X^2+X^3+X^4.
\end{equation}
Since $\Phi_5$ splits into four distinct linear factors over $K$,
\[
 U\cong K^4 *_K K^4.
\]
In $U$, both $a,b$ have fifth power one; moreover
\[
 (X-1)(X^3+2X^2+3X+4)=\Phi_5(X)-5
\]
shows that $a-1,b-1$ are units (as $5\in K^\times$).
Let $L=\ker(U\to U_{\mathrm{ab}})$ be its commutator ideal.
Then $S_0=U/L\cong K^{16}$.  We use $\mathfrak m$ for the
arrow ideal of a completed path algebra.

\begin{lemma}[Formal neighbourhood]\label{lem:quiver}
There is an isomorphism of topological $K$-algebras
\[
 \widehat U_L=\varprojlim U/L^n\ \cong\ \widehat{KQ},
\]
where the vertices of $Q$ are pairs $(i,j)$ with
$1\le i,j\le4$, and one arrow runs from $(i,j)$ to $(k,l)$
exactly when $i\ne k$ and $j\ne l$.  Every vertex has nine
incoming and nine outgoing arrows.
\end{lemma}
\begin{proof}
Write $M=L/L^2$.  The square-zero quotient $U/L^2$ splits
over the separable algebra $S_0$.  Fix such a splitting and
write $e_x$ for the primitive idempotent corresponding to a
character $x$, and $K_x$ for the associated one-dimensional
left $U$-module.  Every extension of $K_x$ by $K_y$ is
annihilated by $L^2$.  The square-zero description gives
\[
 \operatorname{Ext}^1_U(K_x,K_y)
 \cong\operatorname{Hom}_K(e_yMe_x,K).
\]
Indeed, after splitting the extension as an $S_0$-module,
the remaining action is precisely a linear map from
$e_yMe_x$ to $\operatorname{Hom}_K(K_x,K_y)$.
We can therefore compute the corner dimensions of $M$
from upper triangular extensions, modulo changes of splitting.
A character of $S_0$ is specified by a choice of
nontrivial fifth root of unity for each of $a,b$.
For two distinct characters $x=(i,j)$ and $y=(k,l)$, a
square-zero upper triangular extension is determined by the
off-diagonal entries of $a,b$.  The first entry is free exactly
when $i\ne k$, the second exactly when $j\ne l$.
For distinct characters, we quotient out the one-dimensional
space of inner derivations arising from changes of splitting.
For equal characters both entries must
be zero.  Thus the dimensions of the off-diagonal components
of $M$ are one if both labels differ, and zero otherwise;
the diagonal components are zero.  This identifies
$T_{S_0}(M)$ with the ordinary path algebra $KQ$.

For completeness we justify the passage to completions.
A homomorphism $K^r\to E/J$ with nilpotent $J$ lifts to
$E$ by lifting a complete set of orthogonal idempotents:
for a square-zero ideal and $x^2-x=n\in J$ one corrects
$x$ to $x+(1-2x)n$; successive lifting handles higher
nilpotence.  Therefore $U=K^4*_KK^4$ has the lifting
property for nilpotent algebra extensions.  The algebra
$T_{S_0}(M)$ has the same property, since $S_0\cong K^{16}$
and its enveloping algebra $S_0\otimes_K S_0^{\rm op}$ are
semisimple, so the $S_0$-bimodule $M$ is projective.

Put $T=T_{S_0}(M)$, with arrow ideal $\mathfrak m$.
The splitting of $U/L^2$ gives an isomorphism
\[
 f_2:T/\mathfrak m^2\longrightarrow U/L^2
\]
inducing the prescribed identifications in degrees zero and one.
Suppose inductively that a compatible isomorphism
$f_n:T/\mathfrak m^n\to U/L^n$ has been constructed,
where $n\ge2$.  Compose $f_n$ with the quotient map from
the untruncated algebra $T$.  Formal smoothness of $T$
lifts this homomorphism to
$\widetilde f_{n+1}:T\to U/L^{n+1}$.
The image of $\mathfrak m$ lies in $L/L^{n+1}$, because
its reduction modulo $L^n$ does.  Hence this lift annihilates
$\mathfrak m^{n+1}$ and induces a compatible map
\[
 f_{n+1}:T/\mathfrak m^{n+1}\longrightarrow U/L^{n+1}.
\]
It is surjective: it is an isomorphism in degrees zero and
one, and the higher associated graded pieces are generated
by products of degree-one elements.

For the reverse comparison, compose $U\to U/L^n$ with
$f_n^{-1}$.  Formal smoothness of the untruncated algebra $U$
gives a lift $U\to T/\mathfrak m^{n+1}$.
This sends $L$ into $\mathfrak m/\mathfrak m^{n+1}$ and
therefore factors through a map
$U/L^{n+1}\to T/\mathfrak m^{n+1}$, which is surjective
by the same associated graded argument.
Both truncated algebras are finite-dimensional.  For $U/L^{n+1}$,
this follows because $U/L$ is finite-dimensional and $L$ is
generated as a two-sided ideal by $ab-ba$; each layer is a
quotient of a tensor power over $S_0$ of $L/L^2$.
The two surjections force equality of dimensions, so
$f_{n+1}$ is an isomorphism.  The forward lifts are compatible
by construction.  Passing to inverse limits proves
$\widehat U_L\cong\widehat{KQ}$.
\end{proof}

\subsection{Sparse torsion relations}
Enumerate the group words in $a^{\pm1},b^{\pm1}$ as
$w_1,w_2,\ldots$, allowing repetitions.  Their images in
$U/L$ have order dividing five, so $u_n=w_n^5-1\in L$.
Define
\begin{equation}\label{eq:relations}
 r_n=w_n^{5\cdot11^n}-1=(w_n^5-1)^{11^n}
 \in L^{11^n}.
\end{equation}
The equality uses the characteristic-$11$ Frobenius identity
for powers of the \emph{single} element $w_n^5-1$; no
commutativity of $U$ is required.  Let
\[
 \widehat R=\widehat U_L/\overline{(r_1,r_2,\ldots)}.
\]
We use the matrix Golod--Shafarevich--Vinberg inequality for
quotients of completed path algebras, in the following form
proved for \emph{finite} families of (possibly inhomogeneous)
relations in \cite[Theorem 1.2]{ChangTang}:
if the quotient is finite-dimensional, then for every
$0<\rho<1$ there is a vector $H(\rho)$ with all entries
at least one satisfying
\begin{equation}\label{eq:gsv}
 \boldsymbol 1\ \le\
 \bigl(I-\rho M_Q^{\mathsf T}+R(\rho)\bigr)H(\rho).
\end{equation}
Here $M_Q$ is the adjacency matrix and $R(\rho)$ is the
matrix whose $(t,s)$-entry sums $\rho^{d_r}$ over
relations with target $t$, source $s$ and $\mathfrak m$-adic
order $d_r$.

\begin{lemma}[Finite extraction]\label{lem:extract}
Suppose $J$ is the closed two-sided ideal generated by a
possibly countable family $\mathcal R$ in a completed path
algebra $A=\widehat{KQ}$ of a finite quiver.  If $A/J$ is
finite-dimensional, there is a finite subfamily
$\mathcal R_0\subseteq\mathcal R$ such that
$A/\overline{(\mathcal R_0)}$ is finite-dimensional.
\end{lemma}
\begin{proof}
Let $\mathfrak m$ be the arrow ideal.  Finite-dimensionality
implies $\mathfrak m^d\subseteq J$ for some $d$ (the image of
$\mathfrak m$ is a nilpotent ideal in the finite-dimensional
quotient).  The vector space
$\mathfrak m^d/\mathfrak m^{d+1}$ is finite-dimensional.
As $J$ is the closure of the ideal generated by $\mathcal R$,
we can choose finitely many relations giving
\[
 \mathfrak m^d\subseteq I_0+\mathfrak m^{d+1},
 \quad I_0=\overline{(\mathcal R_0)}.
\]
Pass to the complete algebra $A/I_0$ and let $\overline{\mathfrak m}$
be its arrow ideal.  Then
$\overline{\mathfrak m}^{\,d}=
\overline{\mathfrak m}^{\,d+1}$.
The ideal $\overline{\mathfrak m}^{\,d}$ is finitely
generated as a right module by images of the finitely many
paths of length $d$, while $\overline{\mathfrak m}$ lies in
the Jacobson radical (geometric series converge).
Nakayama's lemma gives $\overline{\mathfrak m}^{\,d}=0$.
Thus $A/I_0$ is finite-dimensional.
\end{proof}

\begin{proposition}\label{prop:infinite}
The algebra $\widehat R$ is infinite-dimensional over $K$.
\end{proposition}
\begin{proof}
Suppose it is finite-dimensional.  By Lemma~\ref{lem:extract},
a finite subfamily of the $r_n$, after decomposition into
vertex-to-vertex components, already gives a finite-dimensional
quotient of $\widehat{KQ}$.  Each such component has
$\mathfrak m$-adic order at least $11^n$.
For each fixed source vertex an $r_n$ has at most 16 endpoint
components.  Therefore, for $\rho=1/2$, every column sum of
the relation matrix is bounded by
\[
 16\sum_{n\ge1}2^{-11^n}<\frac1{64}.
\]
Every vertex of $Q$ has indegree and outdegree nine.  Hence
every entry of the row vector
\[
 \boldsymbol1^{\mathsf T}
 \bigl(I-\tfrac12M_Q^{\mathsf T}+R(1/2)\bigr)
\]
is at most $1-9/2+1/64<0$.
This is incompatible with \eqref{eq:gsv} after multiplying
on the left by $\boldsymbol1^{\mathsf T}$, because all
entries of $H(1/2)$ are positive.  Contradiction.
\end{proof}

\subsection{Passing from a completed algebra to automorphisms}
\begin{proof}[Proof of Theorem~\ref{thm:1863}]
Let $B$ be the image of the ordinary algebra $U$ in
$\widehat R$.  This image is dense.  If $B$ were a
finite-dimensional $K$-algebra, it would be a finite set
because $K=\F_{11}$, hence closed in the Hausdorff
$\mathfrak m$-adic topology.  Density would imply
$B=\widehat R$, contrary to Proposition~\ref{prop:infinite}.
Thus $B$ is infinite-dimensional and contains a unit.

The images $a,b\in B$ are units of order exactly five:
$a^5=b^5=1$, and $a-1,b-1$ are units, which precludes
$a=1$ or $b=1$.  Every element of the subgroup
$\Gamma=\langle a,b\rangle\le B^\times$ is represented by
some word $w_n$ and satisfies $w_n^{5\cdot11^n}=1$ by
\eqref{eq:relations}.  Hence $\Gamma$ is periodic.
Since $B$ is generated as a $K$-algebra by the units $a,b$
and their inverses are polynomials in them, we have
$B=K\text{-span}(\Gamma)$.  If $\Gamma$ were finite then
$\dim_K B<\infty$, a contradiction.  Thus $\Gamma$ is infinite.

Take the additive group $A=(B,+)$ and let
$\alpha=L_a$, $\beta=L_b$ be left multiplications.
They are automorphisms of exact order five and
have no nonzero fixed points, since $a-1,b-1$ are
invertible in $B$.  The left regular map
$\Gamma\to\Aut(A)$ is injective (evaluate on $1$), so
$\langle\alpha,\beta\rangle\cong\Gamma$ is infinite
and periodic.
\end{proof}

\section{A commutator counterexample in
\texorpdfstring{$\PSL_2(8)$}{PSL2(8)}}
The condition investigated here is the one stated in the
Kourovka Notebook; compare the related conjectures and
commutator results of Guralnick--Malle and
Guralnick--Tong-Viet--Tracey \cite{GM2013,GTT2025}.

We use $[x,y]=x^{-1}y^{-1}xy$.  As $\F_8$ has even
characteristic, $G=\PSL_2(8)=\SL_2(8)$ is simple.
The polynomial $X^2+X+1$ is irreducible over $\F_8$;
therefore the trace-one matrices of determinant one
constitute a single conjugacy class $C$ of cardinality
$504/9=56$, and each has order three.

\begin{proof}[Proof of Theorem~\ref{thm:1865}]
For any $c,d\in C$, put $s=\tr(cd)\in\F_8$.
The Fricke trace identity
\[
 \tr[c,d]=\tr(c)^2+\tr(d)^2+s^2
 +\tr(c)\tr(d)s-2
\]
gives, in characteristic two,
\begin{equation}\label{eq:fricke}
 t:=\tr[c,d]=s^2+s.
\end{equation}
The additive trace $\operatorname{Tr}_{\F_8/\F_2}(t)$ is zero:
\[
 t+t^2+t^4=s+s^8=0.
\]
Consequently $t=0$ or $t$ is a root of
\begin{equation}\label{eq:tracepoly}
 X^3+X+1=0,
\end{equation}
since the trace-zero subspace of $\F_8$ is precisely the
set of roots of $X(X^3+X+1)$.
If $t=0$, then $s^2+s=0$, so $s\in\{0,1\}$.
Lemma~\ref{lem:tracezero} below proves that $c$ and $d$
commute in these two cases.  Thus $[c,d]=1$.

If $t\ne0$, let $z,z^{-1}$ be the eigenvalues of $[c,d]$
over $\F_{64}$.  Then $z+z^{-1}=t$.  Equation
\eqref{eq:tracepoly} gives
\[
 z^3+z^{-3}=(z+z^{-1})^3+(z+z^{-1})=t^3+t=1.
\]
Thus $z^3$ is a nontrivial cube root of unity and $z$ has
order nine.  Hence $[c,d]$ has order nine.
Finally $O_3(G)=1$ by simplicity, completing the proof.
\end{proof}

\begin{lemma}\label{lem:tracezero}
If $c,d\in\SL_2(8)$ have determinant one, trace one,
and $\tr(cd)\in\{0,1\}$, then $cd=dc$.
\end{lemma}
\begin{proof}
Identify the natural two-dimensional $\F_8$-space with
$E=\F_{64}$ so that $c$ acts by multiplication by a primitive
cube root $\zeta\in E$.  Every $\F_8$-linear endomorphism of
$E$ has a unique expression
\[
 d(x)=ux+vx^8,\qquad u,v\in E.
\]
Its trace and determinant over $\F_8$ are respectively
$u+u^8$ and $uu^8+vv^8$.  In particular
\[
 \tr(d)=u+u^8,\quad
 \tr(cd)=\zeta u+\zeta^8u^8,
 \quad \det(d)=uu^8+vv^8.
\]
The first two equations constitute an invertible
$\F_8$-linear system for $u\in E$: its determinant is
$\zeta+\zeta^8=1$, since $\zeta^8=\zeta^2$.
For $(\tr(d),\tr(cd))=(1,1)$ the unique solution is
$u=\zeta$; for $(1,0)$ it is $u=\zeta^2$.
In either case $uu^8=1$, so $\det(d)=1$ implies
$vv^8=0$, hence $v=0$.  Therefore $d$ is multiplication
by $\zeta$ or $\zeta^2$ and commutes with $c$.
\end{proof}

\begin{remark}[An independent finite certificate]\label{rem:code}
A standard-library-only Python program enumerates all $504$
matrices of $\SL_2(8)$ and all $56^2$ ordered pairs in $C$.
It finds exactly $112$ identity commutators and $3024$
commutators of order nine.  The source is supplied as
\texttt{anc/check\_psl2\_8.py}.  The program checks a finite
instance and is not used to justify any infinite construction.
\end{remark}

\section{Locally nilpotent operator-simple groups: Hall's construction}
We spell out the convention used in Problem 20.120.
A pair $(D,\sigma)$ is \emph{operator-simple} if $D\ne1$ and
\[
 N\trianglelefteq D,\ \sigma(N)=N\quad\Longrightarrow\quad
 N=1\ \text{or}\ N=D.
\]
This is equivalent to $D$ being a minimal nontrivial normal
subgroup of $D\rtimes_\sigma\Z$.

Let $V=\bigoplus_{n\in\Z}\F_p e_n$, and for each $n$ define
a permutation $a_n$ of $V$ by
\[
 a_n(v)=\begin{cases}
 v+e_n,&v_j=0\text{ for all }j>n,\\
 v,&\text{otherwise}.
 \end{cases}
\]
Put $W=\langle a_n:n\in\Z\rangle$ and $D=W'$.
For $u\le v$, the subgroup $\langle a_u,\dots,a_v\rangle$
is the iterated regular wreath product of $v-u+1$ copies
of $C_p$.  It is a finite $p$-group.  Thus $W$ is countable,
locally finite, and locally nilpotent.  It is nonabelian:
$a_0a_1(0)=e_1$ while $a_1a_0(0)=e_0+e_1$.
The shift $s(e_n)=e_{n+1}$ conjugates $a_n$ to $a_{n+1}$;
it induces an automorphism $\sigma$ of $D$.

We now apply precisely the two results of P.~Hall
\cite[Theorems A and D]{Hall1962} needed here:
\begin{enumerate}[label=\textup{(H\arabic*)}]
\item for an ordered wreath power indexed by a set without a
greatest element, its derived group is perfect and every normal
subgroup of that derived group is normal in the whole wreath power;
\item if the order-preserving automorphism group of the indexing
set acts irreducibly in Hall's sense, that derived group is a
minimal nontrivial normal subgroup of the natural split extension.
\end{enumerate}
Here the ordered set is $\Z$, without a greatest element;
its order-preserving automorphism group consists of all
translations, and is irreducible in Hall's sense.  Both
hypotheses therefore hold.  Consequently $D=W'=W''$,
and if $N\trianglelefteq D$ is invariant under $\sigma$,
then (H1) gives $N\trianglelefteq W$, while shift invariance
makes it normal in $W\rtimes\Z$.  By (H2), $N=1$ or $D$.
The specialization to the integer-indexed wreath power and
its extension by the shift is already discussed explicitly
in \cite[Section 1.7]{Hall1962}.

To check the extra fixed-point assertion, define for $g\in W$
its \emph{moving-coordinate set}
\[
 M(g)=\{n\in\Z:\exists v\in V\text{ with }(gv)_n\ne v_n\}.
\]
Every $g\in W$ is a finite word in the $a_n$, so $M(g)$
is finite, and $M(g)=\varnothing$ only if $g=1$.
Moreover $M(\sigma^m(g))=M(g)+m$.  If $\sigma^m(g)=g$
for $m\ne0$, this finite set is invariant under a nonzero
translation and must be empty.  Thus
$C_W(\sigma^m)=1$, hence $C_D(\sigma^m)=1$.
Since $W$ is nonabelian, $D=W'\ne1$; the fixed-point
assertion therefore also implies that $\sigma$ has infinite order.
This proves Theorem~\ref{thm:20120}.

\begin{remark}[Attribution and formulation]
Hall's construction predates the Kourovka Notebook question by
decades.  Theorem~\ref{thm:20120} is an observation about its
operator-group interpretation; it must not be represented as
an independently new group construction.  If Problem 20.120
was intended to use a stronger meaning of ``simple with an
automorphism'', the statement would need to be revisited.
\end{remark}

\section{Verification record and scope}
The short counterexample in Section 4 was checked in two independent
ways: the trace calculation and the exhaustive finite-field
program described in Remark~\ref{rem:code}.  Section 5 uses two
specific theorems of Hall, whose content also appears in a published
biographical account \cite{HallObit}.  The two longer sections rely on
standard results that have been explicitly identified:
Steinberg's tensor-product classification and the relevant
spectrum formulas in \cite{Lytkin2018}, Dickson's subgroup
classification, Dirichlet's theorem, and the completed-quiver
inequality of \cite{ChangTang}.  The complete proofs of these
external theorems are not reproduced here.  The reader should
particularly scrutinize the formal-neighbourhood and
finite-extraction arguments in Section 3 and the use of the
outer field coset in Section 2; these are the main technical
interfaces of the submitted constructions.

\section*{AI assistance disclosure}
An AI language model assisted in preparing the initial manuscript. The submitting author has disclosed the use of generative artificial intelligence and
acknowledges responsibility for the mathematical claims and presentation.

\end{document}